\documentclass[12pt]{amsart}

\usepackage{amsmath,amssymb,amstext,amsthm,bm,amsfonts}
\usepackage{url}
\usepackage{color,graphicx,graphics,xcolor}
\usepackage{hyperref}
\usepackage{mathrsfs}
\usepackage{listings}
\usepackage{appendix}
\usepackage{upref,amsxtra,exscale}
\usepackage{cite}
\usepackage{hyperref}
\usepackage {comment}

\usepackage{chngcntr}
\counterwithin*{equation}{section}

\newtheorem{thm}{Theorem}[section]
\newtheorem{defi}[thm]{Definition}

\newtheorem{cor}[thm]{Corollary}

\newtheorem{rem}[thm]{Remark}

\newtheorem{example}[thm]{Example}

\newcommand{\R}{\mathbb{R}}

\def\qq#1{\qquad \mbox{#1}\quad}
\newcommand{\D}{\displaystyle}

\newcommand{\al}{\alpha}

\newcommand{\e}{\varepsilon}

\newcommand{\ga}{\gamma}

\newcommand{\na}{\nabla}

\newcommand{\Om}{\Omega}
\newcommand{\Omb}{\overline{\Om}}
\newcommand{\p}{\partial}
\newcommand{\s}{\sigma}
\newcommand{\te}{\theta}

\newcommand{\vf}{\varphi}

\title[$L^\infty$ a-priori estimates for $p$-laplacian equations]
{$L^\infty$ a-priori estimates for subcritical $p$-laplacian equations with  a Carathéodory nonlinearity}
\author{Rosa Pardo}
\address[R.~Pardo]
{Departamento de An\'alisis Matem\'atico y Matem\'atica Aplicada,  Universidad Complutense de Madrid, 28040--Madrid, Spain.}
\email{rpardo@ucm.es}
\thanks{The  author is supported by grants  PID2019-103860GB-I00,  MICINN,  Spain, and by UCM-BSCH, Spain, GR58/08, Grupo 920894.}

\date{}

\begin{document}
\maketitle

\begin{abstract}
We present  new $L^\infty$ {\it a priori} estimates for weak solutions of a wide class of subcritical  $p$-laplacian equations in bounded domains.   No hypotheses on the sign of the solutions, neither of the non-linearities are required. This method is based in elliptic regularity for the $p$-laplacian combined either with Gagliardo-Nirenberg or Caffarelli-Kohn-Nirenberg interpolation inequalities.

Let us consider a quasilinear boundary value problem
$ -\Delta_p u= f(x,u),$ in $\Omega,$ with Dirichlet boundary conditions, where $\Omega \subset \R^N $, with $p<N,$ is a bounded smooth domain strictly convex, and $f$  is  a subcritical Carathéodory non-linearity.	
We provide $L^\infty$ a priori estimates for weak  solutions, in terms of their $L^{p^*}$-norm, where $p^*= \frac{Np}{N-p}\ $ is the critical Sobolev exponent. 

By a subcritical non-linearity we mean, for instance,
$|f(x,s)|\le |x|^{-\mu}\, \tilde{f}(s),$
where $\mu\in(0,p),$  and $\tilde{f}(s)/|s|^{p_{\mu}^*-1}\to 0$ as $|s|\to \infty$,  here $p^*_{\mu}:=\frac{p(N-\mu)}{N-p}$ is the critical Sobolev-Hardy exponent.  Our non-linearities includes non-power non-linearities. 

In particular we prove that when $f(x,s)=|x|^{-\mu}\,\frac{|s|^{p^*_{\mu}-2}s}{\big[\log(e+|s|)\big]^\al}\,,$ with $\mu\in[1,p),$  
then, for any $\e>0$ there exists a constant $C_\e>0$ such that for any solution $u\in H^1_0(\Om)$, the following holds
$$
\Big[\log\big(e+\|u\|_{\infty}\big)\Big]^\al\le
C_\e \, \Big(1+\|u\|_{p^*}\Big)^{\, (p^*_{\mu}-2)(1+\e)}\, ,
$$
where  $C_\e$  is independent of the  solution $u$.
\end{abstract}

\medskip

\keywords{A priori estimates, subcritical nonlinearity, changing sign weight, $L^\infty$  a priori bound, singular elliptic equations.} 

\subjclass[2020]{35B45, 35J92, 35A23, 35J25}

\section{Introduction}
Let us consider the following quasilinear boundary value problem involving the $p$-Laplacian
\begin{equation}
\label{eq:ell:pb} 
-\Delta_p u= f(x,u), \quad\mbox{in }\Om, \qquad  u= 0,\quad \mbox{on }\p \Om,
\end{equation}
where $\Delta_p
(u)=$
div$(|Du|^{p-2}Du)$ is the $p$-Laplace operator, $1<p<\infty$, $\Om \subset \R^N $,  $N>p,$ is a bounded,  strictly convex, open subset with $C^{2}$ boundary $\p \Om$, 
and the non-linearity $f:\bar{\Om}\times\R\to \R$ is Carathéodory function (that is, the mapping $f(\cdot,s)$ is measurable for all $s\in\R$,  and the mapping $f(x,\cdot)$ is continuous for almost all $x\in\Om$),  
and {\it subcritical} (see definition \ref{def:sub}). 

\medskip

We analyze the effect of the smoothness of the subcritical non-linearity $f=f(x,u)$ on the $L^\infty(\Omega)$ a priori estimates of  {\it weak solutions} to \eqref{eq:ell:pb}.  
This study is usually focused on positive classical solutions, see the classical references of de Figueiredo-Lions-Nussbaum, and of Gidas-Spruck \cite{deFigueiredo-Lions-Nussbaum, Gidas-Spruck}, see also  \cite{Castro_Pardo_RMC, Castro_Pardo_DCDS}.

\bigskip

A natural question concerning  the class  of uniformly bounded solutions is the following one, 
\begin{enumerate}
\item[(Q1)] {\it  those $L^\infty(\Omega)$ estimates apply also to a bigger class of solutions, in particular to {\it weak solutions} (and to changing sign solutions)?. }
\end{enumerate}
Another natural question with respect to  the class of subcritical non-linearities, can be stated as follows,
\begin{enumerate}
\item[(Q2)] {\it  those $L^\infty(\Omega)$ estimates are valid into a bigger class of non-linearities (not asymptotically powers), and in particular  to non-smooth non-linearities (with possibly changing sign weights)?. }
\end{enumerate}

\medskip

\noindent	In this paper we extend the previous work in \cite{Pardo_arxiv_2022}  for $p=2$, and provide sufficient conditions guarantying uniform $L^\infty(\Omega)$ a priori estimates for any $u\in W_0^{1,p}(\Om)$ weak solution to \eqref{eq:ell:pb}, in terms of their $L^{p^*}(\Omega)$ bounds, in the class of Caratheodory subcritical generalized problems. 
In this class, we  state that any set of weak  solutions uniformly $L^{p^*} (\Om)$ a priori bounded is universally $L^\infty (\Om)$ a priori bounded. Our theorems allow changing sign weights, and singular weights, and also apply to changing sign solutions.

\bigskip

Problem \eqref{eq:ell:pb} with $f(x,s)=|x|^{-\mu} |s|^{q-1}s ,$  $\mu>0$, is known as Hardy's problem, due to its relation with the Hardy-Sobolev inequality.
The Caffarelli-Kohn-Nirenberg interpolation inequality for radial singular weights \cite{Caffarelli-Kohn-Nirenberg}, states that whenever $0\le\mu\le p$, 
\begin{equation}\label{def:2*:mu}
p_{\mu}^*:=\frac{p(N-\mu)}{N-p},
\end{equation} 
is the {\it critical} exponent of the Hardy-Sobolev embedding $\quad
W_0^{1,p}(\Om)\hookrightarrow\ $ $ L^{p_{\mu}^*}(\Om, |x|^{-\mu})\,
$ (this embedding is continuous but not compact). 
For the case $0\le\mu\le p$, using a Pohozaev type identity, Pucci and Servadei \cite{Pucci_Servadei_CPAA_2010} prove some non-existence results in $\R^N$. Some existence and non-existence results for power like nonlinearities can be found in \cite{Arias_Cuesta, Egnell, Ekeland-Ghoussoub, Lan_2014}, see also \cite{Romani} for the case $p=N$.

\bigskip

Usually the term subcritical non-linearity is reserved for power like non-linearities. Next, we expand this concept to nonlinearities including the class $o\big(|s|^{p_{\mu}^*-1}\big)$. 

\begin{defi}
\label{def:sub}
By a {\bf subcritical non-linearity} we mean that $f$ satisfies one of the following two growth conditions:

\begin{enumerate}				
\item[{\rm (H0)}] 
\begin{equation}\label{f:growth}
|f(x,s)|\le |a(x)|\, \tilde{f}(s)
\end{equation}
\end{enumerate}
where $a\in L^r(\Om)$ with $r>N/p$,   $\tilde{f}:\R\to [0,+\infty)$ is continuous and satisfy   
\begin{equation}\label{f:sub}
\tilde{f}(s)>0 \ \text{for}\ |s|>s_0,\ \text{ and }\  \lim_{s \to \pm\infty}  \frac{\tilde{f}(s)}{|s|^{p_{N/r}^*-1}}=0,
\end{equation}
where 
\begin{equation}\label{def:2:r}
p_{N/r}^*:=\frac{p^*}{r'}=p^*\left(1-\frac{1}{r}\right),
\end{equation}
and where $r'$ is the conjugate exponent of $r$,  $1/r+1/r'=1$,
\footnote{Since $r>N/p,$ then  $p_{N/r}^*>p.$
Moreover, thanks to Sobolev embeddings, for any $u\in W_0^{1,p}(\Om)$,
\begin{align}		
\tilde{f}(u)\in L^{\frac{p^*}{p_{N/r}^*-1}}(\Om) & \qq{with} \frac{p_{N/r}^*-1}{p^*}=1-\frac{1}{r}+\frac1N-\frac1p, \nonumber	
\\ 
\text{and}\ f(\cdot,u)\in L^{(p^*)'}(\Om). \nonumber	
\end{align}
}
\ or

\begin{enumerate}				
\item[{\rm (H0)' }] 
\begin{equation}\label{f:growth:rad}
|f(x,s)|\le |x|^{-\mu}\, \tilde{f}(s),
\end{equation}
\end{enumerate}
where $\mu\in(0,p),$  and $\tilde{f}:\R\to[0,+\infty)$   is continuous and satisfy  
\begin{equation}\label{f:sub:rad}
\tilde{f}(s)>0 \ \text{for}\ |s|>s_0,\ \text{ and }\  \lim_{|s| \to \infty} \ \frac{\tilde{f}(s)}{|s|^{p_{\mu}^*-1}}=0
\ \footnote{
Observe that $p_{\mu}^*>p$ for $\mu\in(0,p)$. Let $a(x)=|x|^{-\mu}$, then $a\in L^p(\Om)$ for any $p<N/\mu$. Moreover, the sharp Caffarelli-Kohn-Nirenberg  inequality implies that if $u\in W_0^{1,p}(\Om)$, then
$f(\cdot,u)\in L^{(p^*)'}(\Om)$
}\,.
\end{equation}		
\end{defi}

\medskip

Our analysis shows that non-linearities satisfying either {\rm (H0)}: \eqref{f:growth}-\eqref{f:sub} (either {\rm (H0)' }: \eqref{f:growth:rad}-\eqref{f:sub:rad}),  widen the  class of subcritical non-linearities to non-power non-linearities, sharing with power like non-linearities properties such as $L^\infty$ a priori estimates. 
Our definition of a subcritical non-linearity includes  non-linearities such as 
$$
f^{(1)}(x,s):=
\frac{a(x)|s|^{p_{N/r}^*-2}s}{\big[\log(e+|s|)\big]^\al},
\ \ \text{or}\ \ f^{(2)}(x,s):= 
\frac{|x|^{-\mu}|s|^{p_{\mu}^*-2}s}{\Big[\log\big[e+\log(1+|s|)\big]\Big]^\al},
$$
for any $\al>0$, and  $\mu\in(0,p)$, or any $a\in L^r(\Om),$ with $r>N/p$.

\medskip

In particular, if $f(x,s)=f^{(1)}(x,s)$ with $a\in L^r(\Om)$ for $r\in(N/p,N]$, then for any $\e>0$ there exists a constant $C>0$ depending only on $\e,\ \Om,$ $r$ and $N$ such that for any  $u\in W_0^{1,p}(\Om)$  solution   to \eqref{eq:ell:pb}, the following holds:
\begin{equation*}
\Big[\log\big(e+\|u\|_{\infty}\big)\Big]^\al	\le
C \|a\|_r^{\ 1+\e} \ \Big(1+\|u\|_{p^*}\Big)^{\, (p_{N/r}^*-p)(1+\e)},
\end{equation*}
where  $C$  is independent of the  solution $u$, see Theorem \ref{th:apriori:cnys}. Related results concerning those non-power non-linearities can be found in \cite{Damascelli_Pardo}, and  for $p=2$ in  \cite{Clapp_Pardo_Pistoia_Saldana} analyzing what happen when $\al\to 0$, in \cite{Cuesta_Pardo_MJM_2022} with changing sign weights, in \cite{Mavinga_Pardo17_JMAA} for systems, and in \cite{Pardo_Sanjuan} for the radial case.

\medskip

Moreover, if $f(x,s)=f^{(2)}(x,s)$ 
with $\mu\in [1,p),$ then for any $\e>0$ there exists a constant $C>0$ depending  on $\e$,  $\mu$, $N$,  and $\Om,$ such that for any  $u\in W_0^{1,p}(\Om)$  solution   to \eqref{eq:ell:pb}, the following holds:
\begin{equation*}
\Big[\log\big[e+\log\big(1+\|u\|_{\infty}\big)\big]\Big]^\al\le
C \ 
\Big(1+\|u\|_{p^*}\Big)^{\, (p_{\mu}^*-p)(1+\e)},
\end{equation*}
and where  $C$  is independent of the  solution $u$, see Theorem \ref{th:rad}.

\begin{defi}
\label{def:w:sol}
By a {\it solution} we mean a weak solution $u \in W_0^{1,p}(\Om)$  such that $f(\cdot,u)\in L^{(p^*)'}(\Om),$ and
\begin{equation}\label{weak:sol}
\int_{\Om}|\na|^{p-2}\na u\cdot \na \vf=\int_{\Om} f(x,u)\vf, \quad\quad \forall \vf\in W_0^{1,p}(\Om).
\end{equation}
\end{defi}

\bigskip

To state our main results, for a non-linearity $f$ satisfying
\eqref{f:growth}-\eqref{f:sub}, let us  define
\begin{equation}\label{def:h}
h(s)=h_{N/r}(s):=\frac{|s|^{p_{N/r}^*-1}}{\displaystyle 
	\max\big\{\tilde{f}(-s),\tilde{f}(s)\big\}}\qq{for} |s|> s_0.
\end{equation}
And for a non-linearity $f$ satisfying  \eqref{f:growth:rad}-\eqref{f:sub:rad}, let us now define
\begin{equation}\label{rad:def:h}
h(s)=h_{\mu}(s):=\frac{|s|^{p_{\mu}^*-1}}{\displaystyle 
\max\big\{\tilde{f}(-s),\tilde{f}(s)\big\}}, \qq{for} |s|> s_0.
\end{equation}

By sub-criticallity, (see \eqref{f:sub} or \eqref{f:sub:rad} respectively), 
\begin{equation}\label{h:inf}
h(s)\to\infty \qq{as} s\to\infty.
\end{equation}

\medskip

\noindent Let $u$ be a solution to \eqref{eq:ell:pb}. We estimate $h\big(\|u\|_{\infty}\big)$, in terms of the $L^{p^*}$-norm of $u$.

Our main results are  Theorem \ref{th:apriori:cnys} and Theorem \ref{th:rad},  stated in the following two  subsections for (H0) or (H0)' respectively.

\subsection{Estimates of the $L^\infty$-norm of the solutions to \eqref{eq:ell:pb} in presence of a Caratheodory nonlinearity}\ \\
We assume that the non-linearity 
$f$ satisfies the growth condition {\rm (H0)}, 
and that $\tilde{f}:\R\to (0,+\infty)$  satisfies the following hypothesis:
\begin{enumerate}				
\item[\rm (H1)] there exists a uniform constant $c_0>0$ such that 
\begin{equation}\label{H1}
	\limsup_{s\to +\infty}\ \dfrac{\max_{[-s,s]}\, \tilde{f}}{\max\big\{\tilde{f}(-s),\tilde{f}(s)\big\}}\, \le \, c_0.
	\footnote{Observe that in particular, if  $\tilde{f}(s)$ is monotone, then (H1) is obviously satisfied with $c_0=1$.} 
\end{equation}	
\end{enumerate}

\bigskip

Under hypothesis  {\rm (H0)}-(H1), we establish an estimate for the function $h$ applied to the $L^\infty (\Om)$-norm of any $u\in W_0^{1,p}(\Om)$ 
solution   to \eqref{eq:ell:pb}, in terms of  their $L^{p^*} (\Om)$-norm. 

From now on, $C$ denotes several constants that may change from line to line, and are independent of $u$. 

Our first main results is the following  theorem.

\begin{thm}
\label{th:apriori:cnys}  
Assume that $f:\Omb \times\R\to \R$  is a Carathéodory function  satisfying  {\rm (H0)}-{\rm (H1)}.

Then,      for any  $u\in W_0^{1,p}(\Om)$ weak solution   to \eqref{eq:ell:pb}, the following holds:
\begin{enumerate}				
	\item[{\rm (i)}] either there exists a constant $C>0$ such that $\|u\|_{\infty}\le C$, where  $C$  is independent of the  solution $u$,
	\item[{\rm (ii)}] either 
	for any $\e>0$ there exists a constant $C>0$ such that 
	\begin{equation*}
		h\big(\|u\|_{\infty}\big)
		\le C \|a\|_r^{\ A+\e} \ \Big(1+\|u\|_{p^*}\Big)^{\, (p_{N/r}^*-p)(A+\e)},
	\end{equation*}
	where $h$ is defined by \eqref{def:h}, 
	\begin{equation}\label{def:A}
		A:= 
		\begin{cases}
			1 , &\qq{if}r\le N,\\[.1cm]
			\D \frac{p_{N/r}^*-1}{p_{N/p}^*} ,&\qq{if}r> N,
		\end{cases}	
	\end{equation}
	and $C$ depends only on $\e$, $c_0$ (defined in \eqref{H1}),  $r$, $N$,  and $\Om,$ and it is independent of the  solution $u$.
\end{enumerate}
\end{thm}

As an immediate consequence, as soon as we have a universal {\it a priori} $L^{p^*}$- norm for weak  solutions in $W_0^{1,p}(\Om),$ 
then solutions are a priori universally bounded  in the $L^{\infty}$- norm,  see Corollary \ref{cor:apriori:cnys}.

\subsection{Estimates of the $L^\infty$-norm of the solutions to \eqref{eq:ell:pb} in presence of radial singular weights}\ \\
Now, assuming that $0\in\Omb$ and  that   $|f(x,s)|\le |x|^{-\mu}\, \tilde{f}(s)$ for some $\mu\in (0,p)$, we state our second main result concerning weak solutions for singular subcritical non-linearities, see the following  theorem.

\begin{thm}\label{th:rad}
Assume that $f:\Omb \times\R\to \R$  is a Carathéodory function  satisfying  {\rm (H0)'} and {\rm (H1)}. 

Then,    for any  $u\in W_0^{1,p}(\Om)$ solution   to \eqref{eq:ell:pb}, the following holds:
\begin{enumerate}				
	\item[{\rm (i)}] either there exists a constant $C>0$ such that $\|u\|_{\infty}\le C$, where  $C$  is independent of the  solution $u$,
	\item[{\rm (ii)}] either 
	for any $\e>0$ there exists a constant $C>0$ such that 
	\begin{align*}
		h\big(\|u\|_{\infty}\big)\le C_\e \ \Big(1+\|u\|_{p^*}\Big)^{\, (p_{\mu}^*-p)(B+\e)}\ , 
	\end{align*}
	where $h$ is defined by \eqref{rad:def:h},  
	\begin{equation}\label{def:B}
		B:= 
		\begin{cases}
			\D \frac{p_\mu^*-1}{p_{N/p}^*}
			,&\text{if}\ \mu\in (0,1),\\[.1cm]
			1 ,&\text{if}\ \mu\in [1,p),
		\end{cases}	
	\end{equation}
	and $C$ depends only on $\e$, $c_0$ (defined in \eqref{H1}),  $\mu$, $N$,  and $\Om,$ and it is independent of the  solution $u$.
	\footnote{Observe that $\frac{p_\mu^*-1}{p_{N/p}^*}=1+\frac{p}{p-1}\,\frac{1-\mu}{N}=B$ if $\mu\in(0,1)$.}
\end{enumerate}	
\end{thm}

\bigskip

This results  hold for positive, negative and changing sign non-linearities, and also for positive, negative and changing sign solutions. The techniques and ideas introduced in \cite{Pardo_arxiv_2022} are robust enough to be used for proving analogues of our results in other non-linear problems. Here we present the work for the $p$-Laplacian.
The work for nonlinear boundary conditions is actually in preparation by Chhetri, Mavinga, and the  author.

\bigskip

This paper is organized in the following way. Section \ref{sec:known} collects some well known results.
In  Section \ref{sec:main:I}, using Gagliardo--Nirenberg inequality, we analyze the case when $a\in L^r(\Om)$ with $r>N/p$, see Theorem \ref{th:apriori:cnys}. In  Section \ref{sec:main:II}, we analyze the more involved case of a radial singular weight, see Theorem \ref{th:rad}. It yields on the  Caffarelli-Kohn-Nirenberg inequality.

\section{Preliminaires and known results}
\label{sec:known}

\subsection{Gradient Regularity}\ \\
We are going to use the following result about the summability of the gradient for solutions to equations involving the $p$-Laplace operator. 

\begin{thm}[Gradient Regularity]\label{th:RegularityGradient}
Let $\Omega  $ be a smooth bounded domain in $\R^N$, $N \geq 2$, and let 
$u \in W_0^{1,p}(\Omega ) $, $1 < p < \infty $, be a solution of the problem 
\begin{equation}\label{Cianchi-Iwaniek}
\begin{cases}
-\Delta_p (u) = g &\qquad \text{in }\, \Omega  \\
\qquad\ \  u = \, 0 &\qquad \text{on }\, \partial\Omega, 
\end{cases}
\end{equation}
with $g \in L^q(\Omega )$. 
We assume that 
\begin{equation} 
\begin{cases}
1<q< \infty  & \text{ if } \quad p \geq N, \\
(p^{*})' \le q < \infty &  \text{ if } \quad  1< p < N.
\end{cases}
\end{equation} 
Here 
$p^*= \frac{Np}{N-p}$ is the critical exponent for  Sobolev embedding, and 
$(p^{*})' = \frac{p^*}{p^* -1}= \frac{Np}{Np-N+p}$, is its conjugate exponent.
\begin{itemize}
\item [i)  ] If $ q<N$, then  \ $\ \Vert \na u \Vert_{L^{q^*(p-1)} (\Omega )} \le 
C \Vert g  \Vert_ {L^q (\Omega ) }^{\frac1{(p-1)}}$
\item [ii)  ] If $q\geq N$, then  \ $\ \Vert \na u \Vert_{L^{\sigma  } (\Omega )} \le 
C \Vert g  \Vert_ {L^q (\Omega ) }^{\frac1{(p-1)}}$  \ for any $\sigma < \infty $. 
\end{itemize}
Here $C$ is a constant that depends on $p,N,q$.
\end{thm}

\begin{rem} The exponent $(p^*)'$  is called the duality exponent, and the condition $q \geq (p^{*})'$ if $1<p<N $ guarantees by Sobolev's embeddings that  $g \in L^q(\Omega )$ belongs to the dual space 
$W^{-1, p'}(\Omega )$. If other cases,  we enter into the field of problems with measure data, and  other definitions of solutions have to be considered (see \cite{Boccardo-Gallouet}, \cite{Mingione_2010}).
\end{rem}

The previous theorem follows from different results  proved in several papers (see \cite{Boccardo-Gallouet}, \cite{Byun-Wang-Zhou}, \cite{Cianchi-Mazya_2014}, \cite{DiBenedetto_Nonl-Anal}, \cite{DiBenedetto-Manfredi}, \cite{Iwaniec}, \cite{Mingione_2010}, the survey \cite{Cianchi-Mazya_2015},  and the references therein), where  more general situations are also considered.

\subsection{Improved regularity of the weak solutions}\ \\
We first collect a regularity Lemma for any weak solution to \eqref{eq:ell:pb} with a  non-linearity  of sub-critical growth, in fact weak solutions in $W_0^{1,p}(\Om)$ are in $L^q$  for any finite $q \geq 1$, see 
\cite[Theorem 2.1, Theorem 2.2]{Pucci_Servadei_IUMJ_2008}.

\begin{thm}[Improved regularity]
\label{lem:reg} 
Assume that $u\in W_0^{1,p}(\Om)$ weakly solves \eqref{eq:ell:pb} for a Carathéodory non-linearity $f:\bar{\Om}\times\R\to\R$ with  sub-critical groth, either {\rm (H0)},  either {\rm (H0)'}, see \eqref{f:growth}-\eqref{f:sub} or \eqref{f:growth:rad}-\eqref{f:sub:rad} respectively.

Then, $u\in L^q(\Om)$ for any $1\le q<\infty$.

\smallskip

Moreover, $u\in L^\infty(\Om)$.
\end{thm}
\begin{proof}
We first adapt to the $p$-laplacian the technique used in \cite{deFigueiredo-Lions-Nussbaum} (based in Brezis-Kato, see \cite{Brezis-Kato}) to get the $L^q$ estimates for any finite $q \geq 1$. 

Testing the equation $- \Delta_p u = f(x,u)$ with $|u|^{t}$, $t \geq 1$, we get that 
$$
t \int_{\Omega  } |\na u |^p |u|^{t-1} \, dx = \int  f(x,u) |u|^{t} \,dx . 
$$
Since 
$$
\left|\na\left(|u|^{\frac{p-1+t}{p}}\right)\right|^p 
= \left(\frac{p-1+t}{p} \right)^p |\na u|^p |u|^{t-1},
$$ 
we can write the previous equation as 
\begin{equation}
t \,\left(\frac{p}{(p-1+t)}\right)^p \,\int_{\Omega  }  \left|\na\left(|u|^{\frac{p-1+t}{p}}\right)\right|^p =\int_{\Omega  } f(x,u) |u|^{t}
\end{equation}

\bigskip

(i) We start assuming (H0), see \eqref{f:growth}-\eqref{f:sub}.
By sub-criticallity, (see \eqref{f:sub}), for any $ \e >0 $, there exists $s_{\e}'$ such that 
$$
|f(x,s)|\, |s|^t \le \e\, |a(x)||s|^{p_{N/r}^*-1+t} \qq{if} s \geq s_{\e}',
$$ 
so that denoting by $C_t$ a uniform constant depending also on $t$, we get that 
\begin{align*} 
\int_{\Omega  }  \left|\na\left(|u|^{\frac{p-1+t}{p}}\right)\right|^p 
&\le  C_t \left(C_1+\e \int_{\Omega  } |u|^{p_{N/r}^*-1+t}\, dx \right)\\  
&= C_t  +  \e C_{t} \int_{\Omega  } |u|^{p-1+t} |u|^{p_{N/r}^*-p}  \, dx. 
\end{align*}
By Sobolev's  inequality, and   H\" older's inequality with exponents $\frac{p_{N/r}^*}{p}$, $\frac{p_{N/r}^*}{p_{N/r}^*-p}$, we get that 
\begin{align*} 
&\Big( \int_{\Omega  }   |u|^{\frac{p-1+t}{p} p_{N/r}^* } \, dx \Big)^{\frac{p}{p_{N/r}^*}} 
\le C\int_{\Omega  }  \left|\na\left(|u|^{\frac{p-1+t}{p}}\right)\right|^p \\
&\qquad\qquad \le 
C_t  +  \e C_{t} 
\int_{\Omega  } |u|^{p-1+t} |u|^{p_{N/r}^*-p}  \, dx \\
&\qquad\qquad \le C_t  +  \e C_{t} 
\Big(\int_{\Omega  } |u|^{ \frac{p-1+t}pp_{N/r}^*}\Big)^{\frac{p}{p_{N/r}^*}} \Big(\int_{\Omega  } |u|^{p_{N/r}^*}  \, dx \Big)^{\frac{p_{N/r}^*-p}{p_{N/r}^*}}\\
&\qquad\qquad \le  C_t +  \e C_{t} \Big(\int_{\Omega  } |u|^{ \frac{p-1+t}pp_{N/r}^*}\Big)^{\frac{p}{p_{N/r}^*}}
\end{align*}
Since $u\in W_0^{1,p}(\Om)$, we have that $\int |\na u |^p $ is  bounded. 
Taking $\e $ small we get that 
$\int |u|^{ \frac{p-1+t}pp_{N/r}^*} $ is  bounded for any fixed $1\le t < \infty $, so that   
$\int |u|^{ q} $ is  bounded for any fixed $q \geq p_{N/r}^* $ (and since $\Omega $ is bounded in fact for any 
$q \in [ 1, \infty )$).\\

\bigskip

(ii) We now assume (H0)', see \eqref{f:growth:rad}-\eqref{f:sub:rad}.
By sub-criticallity, (see \eqref{f:sub:rad}), for any $ \e >0 $, there exists $s_{\e}$ such that
$$
|f(x,s)|\, |s|^t \le \e\, |x|^{-\mu}|s|^{p_{\mu}^*-1+t} \qq{if} s \geq s_{\e},
$$ 
so that denoting by $C_t$ a uniform constant depending also on $t$, and   by H\" older's inequality with exponents $\frac{p_{\mu}^*}{p}$, $\frac{p_{\mu}^*}{p_{\mu}^*-p}$, we get that 
\begin{align} \label{grad:ut}
	\nonumber
&\int_{\Omega  }  \left|\na\left(|u|^{\frac{p-1+t}{p}}\right)\right|^p 
\le  C_t \left(C_1+\e \int_{\Omega  }|x|^{-\mu} |u|^{p_{\mu}^*-1+t}\, dx \right)\\  
&\qquad= C_t  +  \e C_{t} \int_{\Omega  } |x|^{-\mu}|u|^{p_{\mu}^*-p} |u|^{p-1+t}  \, dx \nonumber\\  
&\qquad\le C_t  +  \e C_{t} \left(\int_{\Omega  } \left(|x|^{-\mu}|u|^{p_{\mu}^*-p}\right)^{\frac{p^*}{p^*-p}}\right)^{\frac{p^*-p}{p^*}} \left(\int_{\Omega  }|u|^{ \frac{p-1+t}{p} p^*}  \right)^{\frac{p}{p^*}}\nonumber \\ 
&\qquad= C_t  +  \e C_{t} \big\||x|^{-\ga}u\big\|_\rho^{p_\mu^*-p}\,
\left(\int_{\Omega  }|u|^{ \frac{p-1+t}{p} p^*}  \right)^{\frac{p}{p^*}}, 
\end{align}
where $\ga:=\frac{\mu}{p_\mu^*-p}$,  and $\rho:=\frac{(p_\mu^*-p)p^*}{p^*-p}\,$.

Now, since Caffarelli-Kohn-Nirenberg interpolation inequality,
\begin{align} \label{CKN}
\big\||x|^{-\ga}u\big\|_\rho
\le C\, \|\na u\|_{p}^\s\ \|u\|_{p^*}^{1-\s}
\end{align}
where
\begin{equation}\label{coef:CKN}
\frac{1}{\rho}-\frac{\ga}{N}=\s \left(\frac{1}{p}-\frac{1}{N}\right)+ \frac{1-\s}{p^*}=\frac{1}{p^*},
\end{equation}
which trivially holds for any $\s\in(0,1)$. Then, the above can be writen as
\begin{equation}\label{grad:ut:2}
\int_{\Omega  }  \left|\na\left(|u|^{\frac{p-1+t}{p}}\right)\right|^p 
\le  C_t + \e C_{t}\, \big\|\na u\big\|_{p}^{\ p_\mu^*-p}\,
\left(\int_{\Omega  }|u|^{ \frac{p-1+t}{p} p^*}  \right)^{\frac{p}{p^*}}
\end{equation}

Firstly, by Sobolev's  inequality, and secondly by \eqref{grad:ut:2},  we get that 
\begin{align*} 
\Big( \int_{\Omega  }   |u|^{\frac{p-1+t}{p} p^* } \, dx \Big)^{\frac{p}{p^*}} 
& \le C\int_{\Omega  }  \left|\na\left(|u|^{\frac{p-1+t}{p}}\right)\right|^p \\
&\le 
C_t  +  \e C_{t} 
\, \big\|\na u\big\|_{p}^{\ p_\mu^*-p}\,
\left(\int_{\Omega  }|u|^{ \frac{p-1+t}{p} p^*}  \right)^{\frac{p}{p^*}}.
\end{align*}
Since $u\in W_0^{1,p}(\Om)$, we have that $\int |\na u |^p $ is  bounded. 
Taking $\e $ small we get that 
$\int |u|^{ \frac{p-1+t}pp_{\mu}^*} $ is  bounded for any fixed $1\le t < \infty $, so that   
$\int |u|^{ q} $ is  bounded for any fixed $q \geq p_{\mu}^* $ (and since $\Omega $ is bounded in fact for any 
$q \in [ 1, \infty )$).

\bigskip

Finally, combining the above estimates, with the gradient regularity of Theorem \ref{th:RegularityGradient}, and the Sobolev embeddings, we deduce that $u\in L^\infty(\Om)$.
\end{proof}

\section{Carathéodory non-linearities}
\label{sec:main:I}

We start this section with an immediate corollary of Theorem \ref{th:apriori:cnys}:  any sequence of solutions in $W_0^{1,p}(\Om)$, uniformly bounded in the $ L^{p^*} (\Om)$-norm,  is also uniformly bounded in the $ L^{\infty} (\Om)$-norm.

\begin{cor}\label{cor:apriori:cnys}
Let $f:\Omb \times\R\to \R$ be a   Carathéodory function satisfying {\rm (H0)}--{\rm (H1)}. 

Let  $\{u_k\}\subset W_0^{1,p}(\Om)$ be any sequence of solutions to \eqref{eq:ell:pb} such that there exists a constant $C_0>0$ satisfying
\begin{equation*}
	\|u_k\|_{{p^*}} \le C_0.
\end{equation*}

Then, there exists a constant $C>0$ such that
\begin{equation}\label{L:inf:k:0}
	\|u_k\|_{\infty} \le C.
\end{equation} 
\end{cor}

\begin{proof}
We reason by contradiction, assuming that \eqref{L:inf:k:0} does not hold. So, at least for a subsequence again denoted as $u_k$, $\|u_k\|_{\infty} \to\infty$  as $k \to \infty.$
Now part (ii) of the Theorem \ref{th:apriori:cnys} implies that 
\begin{equation}\label{h:bdd:uk:0}
	h\big(\|u_k\|_{\infty}\big)
	\le C.
\end{equation}
From hypothesis {\rm (H0)}(see in particular \eqref{h:inf}), for any $\e>0$ there exists $s_1>0$ such that  $h(s)\ge 1/\e$ for any $s\ge s_1,$ and so $h\big(\|u_k\|_{\infty}\big)\ge 1/\e$ for any $k$ big enough. This contradicts \eqref{h:bdd:uk:0},  ending the proof.	
\end{proof}

\subsection{Proof of Theorem \ref{th:apriori:cnys}}\ \\
\label{sec:proof:apriori:cnys}
The arguments of the proof use Gagliardo-Nirenberg interpolation inequality (see \cite{Nirenberg}), and are inspired in the equivalence between uniform $L^{p^*}(\Om)$ {\it a priori} bounds and uniform $L^\infty (\Om)$ {\it a priori} bounds for solutions to subcritical elliptic  equations, see \cite[Theorem 1.2]{Castro_Mavinga_Pardo} for the quasilinear case and $f=f(u)$, and \cite[Theorem 1.3]{Mavinga_Pardo_MJM} for the $p$-laplacian  and $f=f(x,u)$. We first use elliptic regularity, 
and next, we invoke the Gagliardo-Nirenberg interpolation inequality (see \cite{Nirenberg}).

\begin{proof}
Let   $\{u_k\}\subset W_0^{1,p}(\Om)$ be any sequence of  weak solutions to \eqref{eq:ell:pb}. Since Theorem \ref{lem:reg}, in fact $\{u_k\}\subset W_0^{1,p}(\Om)\cap L^\infty(\Om)$.
If $\|u_k \|_{\infty}\le C,$ then 
(i)  holds.
Now, we  argue on the  contrary, assuming that there exists a sequence $\|u_k \|_{\infty} \to + \infty$ as $k \to \infty.$

We split the proof in two steps. Firstly, we write an $W^{1,q^*(p-1)}$ estimate for $q\in\big(N/p,\min\{r,N\}\big)$,  with $q^*(p-1)>N.$ Secondly, we invoke the Gagliardo-Nirenberg interpolation inequality  for the $L^\infty$-norm in terms of its $W^{1,q^*(p-1)}$-norm and its $L^{p^*}$-norm.

\bigskip

{\it Step 1. $W^{1,q^*(p-1)}$ estimates for $q\in\big(N/p,\min\{r,N\}\big)$.}

\medskip

\noindent Let us denote by 
\begin{equation}\label{def:M:k:f}
M_k:=\max\Big\{\tilde{f}\big(-\|u_k\|_{\infty}\big),\tilde{f}\big(\|u_k\|_{\infty}\big)\Big\}
\ge \frac{1}{2c_0}
\,\max_{[-\|u_k\|_{\infty},\|u_k\|_{\infty}]}\tilde{f},
\end{equation}
where the inequality holds by  hypothesis (H1), see \eqref{H1}.

Let us take $q$ in the interval  $(N/p,N)\cap (N/p,r).$ Growth hypothesis {\rm (H0)}(see \eqref{f:growth}-\eqref{f:sub}), hypothesis (H1) (see \eqref{H1}), and Hölder inequality, yield the following
\begin{align}\label{f:q:k:2}
&\displaystyle \int_{\Om}\left|f\big(x,u_k(x)\big)\right|^q\, dx 
\le\int_{\Om} |a(x)|^q \left(\tilde{f}\big(u_k(x)\big)\right)^{q}\, dx\nonumber\\
&\qquad =\int_{\Om} |a(x)|^q \left(\tilde{f}\big(u_k(x)\big)\right)^{t} \, 
\left(\tilde{f}\big(u_k(x)\big)\right)^{q-t}\, dx\nonumber\\
&\qquad \le  C
\left[\int_{\Om} |a(x)|^q\ \left(\tilde{f}\big(u_k(x)\big)\right)^{t} \,dx\right]
\ M_k^{\ q-t}
\nonumber\\
&\qquad \le  C   
\left(\int_{\Om} |a(x)|^{qs} \,dx\right)^\frac{1}{s} \left(\int_{\Om} \left(\tilde{f}\big(u_k(x)\big)\right)^{ts'} \,dx\right)^\frac{1}{s'}
\ M_k^{\ q-t}
\nonumber\\
&\qquad \le  C\|a\|_r^q\ \Big(\|\tilde{f}(u_k) \|_{\frac{p^*}{p_{N/r}^*-1}}\Big)^{t}\
M_k^{\ q-t},
\end{align}
where $\frac{1}{s}+\frac{1}{s'}=1$, $qs=r$, $C=c_0^{q-t}$  (for $c_0$ defined in \eqref{H1}), and $ts'=\frac{p^*}{p_{N/r}^*-1}$, so 
\begin{align}\label{def:t}
t&:=\frac{p^*}{p_{N/r}^*-1}\left(1-\frac{q}{r}\right)<q \\
&\iff \frac1{q}-\frac1{r}<\frac{p_{N/r}^*-1}{p^*}
=1-\frac1{r}-\frac1p+\frac1N \nonumber\\
&\iff \frac1{q}<1-\frac1p+\frac1N\iff \frac1{q}<1-\frac1{p^*}=\frac1{(p^*)'}\ \checkmark \nonumber
\end{align}
since $p/N<1-\frac1{p^*}\iff p<N \checkmark$, and   $q>N/p>(p^*)'.$

By the gradient regularity for the $p$-laplacian (see Theorem \ref{th:RegularityGradient}) we have that
\begin{equation}\label{ell:reg:k}
\|\na u_k\|_{L^{q^*(p-1)} (\Omega)}\le     C\
\left\|f\big(\cdot,u_k(\cdot)\big)\right\|_q^{\ \frac1{p-1}},
\end{equation}
where $1/q^*=1/q-1/N$, and $C=C(c_0,N,p,q,|\Om|)$ and it is independent of $u.$ 
Since $q>N/p$, then
\begin{equation}\label{def:r:k}
r:=q^*(p-1)>N.
\end{equation}

Now, substituting \eqref{f:q:k:2} into  \eqref{ell:reg:k}
\begin{equation*}
\|\na u_k\|_{L^{q^*(p-1)} (\Omega)}
\le   C\
\left(\|a\|_r\ \Big(\|\tilde{f}(u_k) \|_{\frac{p^*}{p_{N/r}^*-1}}\Big)^\frac{t}{q}\
M_k^{\ 1-\frac{t}{q}}\right)^{\ \frac1{p-1}},
\end{equation*}
Observe that since $q>N/p$, then $q^*(p-1)>N.$

\bigskip

{\it Step 2. Gagliardo-Nirenberg interpolation inequality.}

\medskip

\noindent Thanks to the Gagliardo-Nirenberg interpolation inequality, there exists a constant $C=C(N,q,|\Om|)$ such that 

\begin{equation*}
\|u_k\|_{\infty}\le C \|\na u_k\|_{q^*(p-1)}^\s \ \|u_k\|_{p^*}^{1-\s}
\end{equation*}
where
\begin{align}\label{Gag-Nir:2}
\nonumber\frac{1-\s}{p^*}
&=\s \left(\frac1{N}-\frac{1}{q^*(p-1)}\right)\\
&=\frac{\s}{p-1}\left(\frac{p-1}{N}-\frac{1}{q}+\frac{1}{N}\right)
=\frac{\s}{p-1}\left(\frac{p}{N}-\frac{1}{q}\right)\nonumber\\
&=\frac{\s}{p-1}\left[1-\frac{1}{q}-p\left(\frac{1}{p}-\frac{1}{N}\right)\right]=\frac{\s}{(p-1)p^*}\left(p_{N/q}^*-p\right).
\end{align}
Hence
\begin{equation}\label{Gag-Nir:3}
\|u_k\|_{\infty}\le C 
\left[
\|a\|_r\ \Big(\|\tilde{f}(u_k) \|_{\frac{p^*}{p_{N/r}^*-1}}\Big)^{\frac{t}{q}}\
M_k^{\ 1-\frac{t}{q}}	
\right]^\frac{\s}{p-1}
\ \|u_k\|_{p^*}^{1-\s},
\end{equation}
where $C=C(c_0,r,N,q,|\Om|)$.

From definition of $M_k$ (see \eqref{def:M:k:f}), and definition of $h$ (see \eqref{def:h}),
we deduce that
\begin{equation*}
M_k	=\frac{\|u_k\|_{\infty}^{\ p_{N/r}^*-1}}{h\big(\|u_k\|_{\infty}\big)}.
\end{equation*}
From \eqref{Gag-Nir:2}
\begin{align}\label{Gag-Nir:4}
\nonumber\frac{1}{\s}
&=1+p^*\left(\frac1{N}-\frac{1}{q^*(p-1)}\right)
\\ \nonumber
&=\frac{1}{(N-p)q(p-1)}\left[(N-p)q(p-1)+pq(p-1)-p(N-q)\right]
\\ \nonumber
&=\frac{1}{(N-p)q(p-1)}\left[Nq(p-1)-p(N-q)\right]
\\ \nonumber
&=\frac{1}{(N-p)q(p-1)}\left[Np(q-1)-q(N-p)\right]
\\ 
&=\frac{1}{p-1}\left[p^*-\frac{p^*}{q}-1\right]
=\frac{1}{p-1}\left(p_{N/q}^*-1\right).
\end{align}
Moreover, since definition of  $t$ (see  \eqref{def:t}),
and definition of $p_{N/r}^*$  (see  \eqref{def:2:r}
\begin{align}
\label{def:be:0}
1-\frac{t}{q}&
=\frac{p^*\left(1-\frac{1}{r}\right)-1-p^*\left(\frac1{q}-\frac{1}{r}\right)}{p_{N/r}^*-1}
=\frac{p_{N/q}^*-1}{p_{N/r}^*-1},	
\end{align}
which, joint with \eqref{Gag-Nir:4}, yield
\begin{equation*}\label{def:be:2}
\frac{\s}{p-1}\left[1-\frac{t}{q}\right]{(p_{N/r}^*-1)}=1.
\end{equation*}
Now \eqref{Gag-Nir:3} can be rewritten as
\begin{equation*}
h\big(\|u_k\|_{\infty}\big)^{\ (1-\frac{t}{q})\frac{\s}{p-1}}\le C 
\left[
\|a\|_r\ \Big(\|\tilde{f}(u_k) \|_{\frac{p^*}{p_{N/r}^*-1}}\Big)^{\frac{t}{q}}\ \right]^\frac{\s}{p-1}
\ \|u_k\|_{p^*}^{1-\s},
\end{equation*}
or equivalently
\begin{equation*}
h\big(\|u_k\|_{\infty}\big)\le C 
\|a\|_r^{\ \te} \ \Big(\|\tilde{f}(u_k) \|_{\frac{p^*}{p_{N/r}^*-1}}\Big)^{\te-1}\		
\ \|u_k\|_{p^*}^{\ \vartheta},
\end{equation*}
where 
\begin{align}
\label{def:al}
\te &:=(1-t/q)^{-1}=\frac{p_{N/r}^*-1}{p_{N/q}^*-1},\\
\label{def:al:be}
\vartheta &:=\frac{1-\s}{\s}(1-t/q)^{-1}(p-1)
=  \te\ (p_{N/q}^*-p),
\end{align}
see \eqref{def:be:0} and \eqref{Gag-Nir:2}.
Observe that since $q< r$, then $\te> 1$. Moreover, since  \eqref{def:al}
\begin{equation}\label{al:be}
\te-1 =\frac{p_{N/r}^*-p_{N/q}^*}{p_{N/q}^*-1}.	
\end{equation}

Furthermore, from sub-criticallity, see \eqref{f:sub}
$$
\int_\Om |\tilde{f}(u_k) |^{\frac{p^*}{p_{N/r}^*-1}}
\le C \left(1+\int_\Om |u_k|^{{p^*}}\, dx\right),
$$
so
\begin{equation*}
\|\tilde{f}(u_k) \|_{\frac{p^*}{p_{N/r}^*-1}}
\le C
\left(1+\|u_k\|_{{p^*}}^{p_{N/r}^*-1}\right).
\end{equation*}
Consequently
\begin{equation*}
h\big(\|u_k\|_{\infty}\big)\le C 
\|a\|_r^{\ \te} \ \Big(1+\|u_k\|_{{p^*}}^{\Theta}\Big),
\end{equation*}
with
$$
\Theta:=(p_{N/r}^*-1)(\te-1)+\vartheta=(p_{N/r}^*-p)\te,
$$
where we have used \eqref{al:be}, \eqref{def:al:be}, and \eqref{def:al}.

Fixed $N>p$ and $r>N/p$, 
the function $q\to \theta=\theta(q)$ for $q\in \big(N/p,\min\{r,N\}\big)$, is decreasing, so 
$$
\inf_{q\in (N/p,\min\{r,N\})} \theta(q)=\theta\big(\min\{r,N\}\big)=A:=\begin{cases}
1,&\text{if}\ r\le N,\\
\frac{p_{N/r}^*-1}{p_{N/p}^*},&\text{if}\ r> N.
\end{cases}
$$
Finally, and since the infimum is not attained in $\big(N/p,\min\{r,N\}\big)$, 
for any $\e>0$,
there exists a constant $C>0$ such that
\begin{equation*}
h\big(\|u_k\|_{\infty}\big)\le 
C\ \|a\|_r^{\ A+\e} \ \Big(1+\|u_k\|_{p^*}^{\, (p_{N/r}^*-p)(A+\e)}\Big),
\end{equation*}
where  $C=C(\e,c_0,r,N,|\Om|),$ ending the proof.
\end{proof}

\section{Radial singular weights}
\label{sec:main:II}

We start this section with their corresponding immediate corollary of Theorem \ref{th:rad}:  any sequence of solutions in $W_0^{1,p}(\Om)$, uniformly bounded in the $ L^{p^*} (\Om)$-norm,  is also uniformly bounded in the $ L^{\infty} (\Om)$-norm. Theire proof is identical to that of Corollary \ref{cor:apriori:cnys}, we omit it.

\begin{cor}\label{cor:apriori:rad}
Let $f:\Omb \times\R\to \R$ be a   Carathéodory function satisfying {\rm (H0)'}--{\rm (H1)}. 

Let  $\{u_k\}\subset W_0^{1,p}(\Om)$ be any sequence of solutions to \eqref{eq:ell:pb} such that there exists a constant $C_0>0$ satisfying
\begin{equation*}
	\|u_k\|_{{p^*}} \le C_0.
\end{equation*}

Then, there exists a constant $C>0$ such that
\begin{equation}\label{L:inf:k:rad}
	\|u_k\|_{\infty} \le C.
\end{equation} 
\end{cor}

\subsection{Proof of  Theorem \ref{th:rad}}\ \\
\label{sec:proof:rad}
We split the proof in two steps, in the first one we use elliptic regularity, 
in the second one, the Caffarelli-Kohn-Nirenberg interpolation inequality for singular weights (see \cite{Caffarelli-Kohn-Nirenberg}).
\begin{proof}
Let   $\{u_k\}\subset W_0^{1,p}(\Om)$ be any sequence of  solutions to \eqref{eq:ell:pb}.  Since Theorem \ref{lem:reg}, $\{u_k\}\subset W_0^{1,p}(\Om)\cap L^\infty(\Om)$.
If $\|u_k \|_{\infty}\le C,$ then 
(i)  holds.

Now, we  argue on the  contrary, assuming that there exists a sequence $\{u_k\}\subset W_0^{1,p}(\Om)$ of  solutions to \eqref{eq:ell:pb}, such that $\|u_k \|_{\infty} \to + \infty$ as $k \to \infty.$ 

By Morrey's Theorem (see \cite[Theorem 9.12]{Brezis}), observe that also 
\begin{equation}\label{Morrey:th}
\|\na u_k \|_{q} \to + \infty \qq{as} k \to \infty,
\end{equation}
for any $q>N$.

\bigskip

{\it Step 1. $W^{1,q^*(p-1)}$ estimates for $q\in\big(N/p,\min\{N,N/\mu\}\big)$.}

\medskip

\noindent As in the proof of Theorem \eqref{th:apriori:cnys}, let us denote by 
\begin{equation}\label{def:M:k:f:rad}
M_k:=\max \Big\{\tilde{f}\big(-\|u_k\|_{\infty}\big),\tilde{f}\big(\|u_k\|_{\infty}\big)\Big\}\ge \frac{1}{2c_0}
\,\max_{[-\|u_k\|_{\infty},\|u_k\|_{\infty}]}\tilde{f},
\end{equation}
where the inequality is due to hypothesis (H1), see \eqref{H1}.

Let us take $q$ in the interval  $(N/p,N)\cap (N/p,N/\mu).$ Using growth hypothesis {\rm (H0)' }(see \eqref{f:growth:rad}), hypothesis (H1) (see \eqref{H1}), and Hölder inequality, we deduce
\begin{align*}
&\displaystyle \int_{\Om}\left|f\big(x,u_k(x)\big)\right|^q\, dx 
\le\int_{\Om} |x|^{-\mu q} \left(\tilde{f}\big(u_k(x)\big)\right)^{q}\, dx\nonumber\\
&\qquad =\int_{\Om} |x|^{-\mu q} \left(\tilde{f}\big(u_k(x)\big)\right)^\frac{t}{p_{\mu}^*-1} \, 
\left(\tilde{f}\big(u_k(x)\big)\right)^{q-\frac{t}{p_{\mu}^*-1}}\, dx\nonumber\\
&\qquad \le  C\ 
\left[\int_{\Om} |x|^{-\mu q}\  \big(1+u_k(x)^{t}\big) \,dx\right]
\ M_k^{\ q-\frac{t}{p_{\mu}^*-1}}
\nonumber\\
&\qquad \le C\,\Big(1+\big|\,|x|^{-\ga}\ u_k\,\big|_{t}^{\ t}\Big)\ \
M_k^{\ q-\frac{t}{p_{\mu}^*-1}},
\end{align*}
where $\ga=\frac{\mu q}{t}$,  $t\in \big(0, q\big(p_{\mu}^*-1\big)\big)$, $C=c_0^{q-\frac{t}{p_{\mu}^*-1}}$  (for $c_0$ defined in \eqref{H1}), and where $M_k$ is defined by \eqref{def:M:k:f:rad}.

Since elliptic regularity see Theorem \ref{th:RegularityGradient},   we have that
\begin{equation}\label{rad:ell:reg}
\|\na u_k\|_{q^*(p-1)}
\le   C\
\left[\Big(1+\big|\,|x|^{-\ga}\ u_k\,\big|_{t}^{\ t}\Big)^{\,\frac{1}{q}}\
M_k^{\ 1-\frac{t}{q(p_{\mu}^*-1)}}\right]^\frac1{p-1},
\end{equation}
where $1/q^*=1/q-1/N$ (since $q>N/p$, then $q^*(p-1)>N$), and $C=C(N,q,|\Om|).$

\bigskip

{\it Step 2. Caffarelli-Kohn-Nirenberg interpolation inequality.}

\medskip

\noindent Since the Caffarelli-Kohn-Nirenberg interpolation inequality, 
there exists a constant $C>0$ depending on the parameters $N,\ q,\ \mu,$ and $t$, such that 
\begin{equation}\label{rad:CKN}
\big|\,|x|^{-\ga}\ u_k\,\big|_{t}
\le C \|\na u_k\|_{q^*(p-1)}^\te \ \|u_k\|_{p^*}^{1-\te},
\end{equation}
where
\begin{align}\label{rad:CKN:2}
\nonumber\frac{1}{t}- \frac{\mu q}{N t}
&= -\te \left( \frac1{N}-\frac{1}{q^*(p-1)}\right)+(1-\te)\frac{1}{p^*}\\
& =\frac{1}{p^*}-\te \left( \frac1{p}-\frac{1}{q^*(p-1)}\right)
\nonumber \\
&=\frac{1}{p^*}- \frac{\te}{p-1}\left(1-\frac1{p} -\frac{1}{q}+\frac1{N}\right)=\frac{1}{p^*}- \frac{\te(p_{N/q}^*-1)}{(p-1)p^*}
,
\end{align}
for $\te\in \left(0,\frac{p-1}{p_{\mu}^*-1}\right)$, and $\frac{\mu q}{t}=\ga$. 
\footnote{Observe that since  $t<q(p_{\mu}^*-1)$, then the r.h.s. of \eqref{rad:CKN:2} is bounded from below, 
$$	
\frac{1}{p^*}- \frac{\te}{(p-1)p^*}(p_{N/q}^*-1)
>\frac{1}{p_{\mu}^*-1}\left(\frac{1}{q}-\frac{\mu}{N}\right),
$$
so 
$$	
\frac{\te}{(p-1)}(p_{N/q}^*-1)<1-\frac{1}{1-\frac{\mu}{N}}\left(\frac{1}{q}-\frac{\mu}{N}\right)=\frac{1-\frac{1}{q}}{1-\frac{\mu}{N}},
$$
and we get the upper bound of $\te$. 
}

Substituting now  \eqref{rad:CKN} into \eqref{rad:ell:reg} we can write
\begin{equation*}
\|\na u_k\|_{q^*(p-1)}\le C \left[
\Big(1+\|\na u_k\|_{q^*(p-1)}^{\te t}\ \|u_k\|_{p^*}^{(1-\te)t}\Big)^\frac{1}{q}\
M_k^{\ 1-\frac{t}{q(p_{\mu}^*-1)}}\right]^\frac1{p-1}.
\end{equation*}
Now, dividing by $\|\na u_k\|_{q^*(p-1)}^{\ \te t/q(p-1)}$ and using \eqref{Morrey:th} we obtain
\begin{equation}\label{rad:CKN:4}
\Big(\|\na u_k\|_{q^*(p-1)}\Big)^{1-\frac{\te t}{q(p-1)}}\le C\left[
\Big(1+\|u_k\|_{p^*}^{\, \frac{(1-\te)t}{q}}\Big)\
M_k^{\, 1-\frac{t}{q(p_{\mu}^*-1)}}\right]^\frac1{p-1}.
\end{equation}

Let us check that
\begin{equation}\label{rad:CKN:4:1}
1-\frac{\te\, t}{q(p-1)}>0 \qq{for any}t<q\big(p_{\mu}^*-1\big).
\end{equation}
Indeed, observe first that \eqref{rad:CKN:2} is equivalent to 
\begin{equation}\label{rad:CKN:2:1}
\te=\frac{\frac{1}{p^*} -\frac{1}{t}+\frac{\mu q}{N t}}{\frac1{p}-\frac{1}{q^*(p-1)}}\ .
\end{equation}
Moreover, from \eqref{rad:CKN:2:1}
\begin{equation}\label{rad:CKN:4:2}
\frac{\te\, t}{q(p-1)}=\frac{\frac{1}{q}\left(\frac{t}{p^*} -1\right)+ \frac{\mu }{N}}{\frac{p-1}{p}-\frac{1}{q^*}}=\frac{\frac{1}{q}\left(\frac{t}{p^*} -1\right)+ \frac{\mu }{N}}{1-\frac{1}{q}-\frac{1}{p^*}}\ ,
\end{equation}
consequently
\begin{align*}
\frac{\te\, t}{q(p-1)}<1 &\iff  \frac{1}{q}\left(\frac{t}{p^*} -1\right)+ \frac{\mu }{N}< 1-\frac{1}{q}-\frac{1}{p^*}\\
&\iff
\frac{1}{q}\frac{t}{p^*}< 
1-\frac1{p^*}-\frac{\mu }{N}\\
&\iff
\frac{t}{q}<p^*\left(1- \frac{\mu }{N }\right)-1=p_{\mu}^*-1
\\
&\iff
t<q\big(p_{\mu}^*-1\big),
\end{align*}
so, \eqref{rad:CKN:4:1} holds.

Consequently,
\begin{equation}\label{rad:CKN:5}
\|\na u_k\|_{q^*(p-1)}\le C
\Big(1+ \|u_k\|_{p^*}^{A_0}\Big)\
M_k^{B_0},
\end{equation}
where
\begin{equation}\label{def:A:B}
A_0:=\frac{\frac{(1-\te)t}{q(p-1)}}{1-\frac{\te\, t}{q(p-1)}}
,\qquad 
B_0:=\frac{\left(1-\frac{t}{q(p_{\mu}^*-1)}\right)\frac{1}{p-1}}{1-\frac{\te\, t}{q(p-1)}}.
\end{equation} 

\bigskip

{\it Step 3. Gagliardo-Nirenberg interpolation inequality.}

\medskip

\noindent Thanks to the Gagliardo-Nirenberg interpolation inequality 
(see \cite{Nirenberg}), there exists a constant $C=C(N,q,|\Om|)$ such that 

\begin{equation}\label{rad:Gag-Nir}
\|u_k\|_{\infty}
\le C \|\na u_k\|_{q^*(p-1)}^\s \ \|u_k\|_{p^*}^{1-\s},
\end{equation}
where
\begin{align}\label{rad:Gag-Nir:2}
\frac{1-\s}{p^*}&=\s \left[\frac1{N}-\frac{1}{(p-1)q^*}\right].
\end{align}
Hence, substituting \eqref{rad:CKN:5} into \eqref{rad:Gag-Nir} we deduce
\begin{equation}\label{rad:Gag-Nir:12}
\|u_k\|_{\infty}\le C 
\ \Big(1+ \|u_k\|_{p^*}^{\ \s\,A_0+1-\s}\Big)
\ M_k^{\ \s B_0}		
.
\end{equation}

\medskip

From definition of $M_k$ (see \eqref{def:M:k:f}), and 
of $h$  (see \eqref{rad:def:h}),
we obtain
\begin{equation}\label{rad:def:M:k:f:10}
M_k	=\frac{\|u_k\|_{\infty}^{\ p_{\mu}^*-1}}{h\big(\|u_k\|_{\infty}\big)}.
\end{equation}
Now we  chech that
\begin{equation}\label{sBp*}
\s B_0\, (p_{\mu}^*-1)=1.
\end{equation}
Indeed, from \eqref{rad:Gag-Nir:2}
\begin{align}\label{rad:Gag-Nir:13}
\frac{1}{\s}&=1+p^*\left(\frac1{N}-\frac{1}{q^*(p-1)}\right)
=\frac{p^*}{p}-\frac{p^*}{q^*(p-1)}\nonumber\\
&= \frac{p^*}{p-1}\left[1-\frac{1}{p}-\frac{1}{q}+\frac{1}{N}\right]=\frac{1}{p-1}\big(p_{N/q}^*-1\big).
\end{align}
From \eqref{rad:CKN:4:2}, 
we deduce
\begin{align}\label{rad:CKN:4:3}
\nonumber 1-\frac{\te\, t}{q(p-1)} &=\frac{1-\frac{1}{q}-\frac{1}{p^*}-\frac{1}{q}\left(\frac{t}{p^*} -1\right)- \frac{\mu }{N}}{1-\frac{1}{q}-\frac{1}{p^*}}\\
&=\frac{\left(1-\frac{\mu }{N}\right)-\frac{1}{p^*}-\frac{t}{qp^*} }{1-\frac{1}{q}-\frac{1}{p^*}}=\frac{p_{\mu}^*-1-\frac{t}{q}}{p^*\left(1-\frac{1}{q}\right)-1}
=\frac{p_{\mu}^*-1-\frac{t}{q}}{p_{N/q}^*-1}\ .
\end{align}
Moreover, since  \eqref{rad:CKN:4:3},
\begin{align}\label{rad:def:be}
&\left(1-\frac{t}{q(p_{\mu}^*-1)}\right){\big(p_{\mu}^*-1\big)}
\frac1{\Big(1-\frac{\te\, t}{q(p-1)}\Big)}
\\ 
&\qquad =\left( p_{\mu}^*-1-\frac{t}{q}\right)
\frac1{\Big(1-\frac{\te\, t}{q(p-1)}\Big)}=p_{N/q}^*-1.
\nonumber	
\end{align}
Hence
\begin{equation}\label{Bp*}
B_0\, (p_{\mu}^*-1)=\frac{p_{N/q}^*-1}{p-1}.
\end{equation}
Taking into account \eqref{rad:Gag-Nir:13} and \eqref{Bp*},  we deduce that \eqref{sBp*} holds.

\smallskip

Consequently, we can rewrite \eqref{rad:Gag-Nir:12} in the following way
\begin{equation}\label{rad:Gag-Nir:14}
h\big(\|u_k\|_{\infty}\big)^{\ \frac1{p_{\mu}^*-1}}\le C \ \Big(1+ 
\|u_k\|_{p^*}^{\ \s\,A_0+1-\s}\Big),
\end{equation}
or equivalently
\begin{equation}\label{rad:Gag-Nir:15}
h\big(\|u_k\|_{\infty}\big)\le C\ \Big(1+  \|u_k\|_{p^*}^{\ \Theta}\Big),
\end{equation}
where
\begin{align*}
\Theta &:=
\big(p_{\mu}^*-1\big)\left[1+\s\,\frac{\frac{t}{q(p-1)}-1}{1-\frac{\te\, t}{q(p-1)}}\right].	
\end{align*}

\medskip

Since \eqref{rad:Gag-Nir:13}-\eqref{rad:CKN:4:3}, $\s\left(1-\frac{\te\, t}{q(p-1)}\right)^{-1}=(p-1)\, (p_{\mu}^*-1 -\frac{t}{q})^{-1}$, and substituting it into the above equation we obtain
\begin{align*}
\Theta &=\big(p_{\mu}^*-1\big)\left( \frac{p_{\mu}^*-p}{p_{\mu}^*-1 -\frac{t}{q}}\right).
\end{align*}
Fixed $p>N$ and $\mu\in (0,p)$, the function $(t,q)\to \Theta=\Theta(t,q)$ for $(t,q)\in \big(0, q(p_{\mu}^*-1)\big)\times \big(N/p,\min\{N,N/\mu\}\big)$, is increasing in $t$ and decreasing in $q$.

For $\mu\in[1,p)$,  $\min\{N,N/\mu\}=N/\mu$. Equation \eqref{rad:CKN:2} with $q=q_k=N/\mu(1-1/k)\to N/\mu$, $t=t_k=\frac{1}{2kp^*}$  and  $\te=\te_k=\frac{p-1}{2(p_{N/q_k}^*-1)}$ is satisfied.
Hence, when $\mu\in[1,p)$,
$$
p_{\mu}^*-p\le \underset{t\in\big(0, (p_{\mu}^*-1)q\big),\, q\in 
\big(\frac{N}{p},\frac{N}{\mu}\big)}{\inf}\ \Theta(t,q)
\le\Theta\big(t_k,q_k\big)\to p_{\mu}^*-p.
$$

On the other hand, for $\mu\in(0,1)$,  $\min\{N,N/\mu\}=N$. For any $\e_k\to 0$,  equation \eqref{rad:CKN:2} with $q=q_k=N(1-\e_k)\to N$, and $t=t_k\to t_0\in [0,(p_\mu^*-1)N]$,  yields   $\te=\te_k=\frac{(p-1)p^*}{p_{N/q_k}^*-1}\left[\frac{1}{p^*}-\frac{1}{t_k}\big(1-\mu(1-1/k)\big)\right]\to p\left[\frac{1}{p^*}-\frac{1}{t_0}(1-\mu)\right]\ge 0$, so $t_0\ge p^*(1-\mu)$. Hence, when $\mu\in(0,1)$,
$$
\underset{t\in[p^*(1-\mu), (p_{\mu}^*-1)N),\, 
q\in \big(\frac{N}{p},N\big)}{\inf}\,  \Theta(t,q)
=\Theta(p^*(1-\mu),N)=(p_{\mu}^*-p)B,
$$
where $B$ is defined by \eqref{def:B}.

Since the infimum is not attained, for any $\e>0$,
there exists a constant $C=C(\e,c_0,\mu,N,\Om)$ such that 
\begin{equation}\label{rad:Gag-Nir:16}
h\big(\|u_k\|_{\infty}\big)\le C \ 
\Big(1+\|u_k\|_{p^*}^{\, (p_{\mu}^*-p)(B+\e)}\Big) ,
\end{equation}
which ends the proof.
\end{proof}

\def\cprime{$'$}

\end{document}